\documentclass[a4paper,10pt]{nyjm}
\usepackage{amsmath}
\usepackage{amsthm}
\usepackage{graphicx}
\usepackage{float}
\title{Further notes on a family of continuous, non-differentiable functions}
\author{James McCollum}
\address{Division of Mathematics and Computer Science, University of Maine at Farmington, Farmington, ME 04938}
\email{james.mccollum@maine.edu}
\date{July 2011}
\thanks{This work was partially supported by a Michael D. Wilson Scholarship from the University of Maine at Farmington.}
\begin{document}

\begin{abstract}
We examine a parametrized family of functions $F_a$, all of which are continuous and some of which are nowhere or almost nowhere differentiable, we explore the behavior of $F_a'$ and $F_a''$ almost everywhere for different values of $a$, focusing on specific questions regarding $F_a$'s differentiability for certain $a$, and we calculate the Hausdorff dimension of the graphs of all $F_a$.
\end{abstract}

\maketitle

\section{Introduction}
In this paper, we will examine a family of functions $F_a$ that are continuous for all $a\in(0,1)$, but nowhere or almost nowhere differentiable for certain $a$.  Their construction, described by Okamoto \cite{Okamoto05} for all $a$,  is based largely on the one given by Bourbaki \cite{Bourbaki04} for $a=2/3$.

Let $F_a$ be defined inductively over $[0,1]$ by iterations $f_i$ for $i\geq0$ as follows: $f_0(x)=x$ for all $x \in [0,1]$, every $f_i$ is continuous on $[0,1]$, every $f_i$ is affine in each subinterval $[k/3^i,(k+1)/3^i]$ where $k\in\{0,1,2,\ldots,3^i-1\}$, and
\begin{flalign}
&f_{i+1}\left(\frac{k}{3^i}\right)=f_i\left(\frac{k}{3^i}\right),\\
&f_{i+1}\left(\frac{3k+1}{3^{i+1}}\right)=f_i\left(\frac{k}{3^i}\right)+a\left[f_i\left(\frac{k+1}{3^i}\right)-f_i\left(\frac{k}{3^i}\right)\right],\\
&f_{i+1}\left(\frac{3k+2}{3^{i+1}}\right)=f_i\left(\frac{k}{3^i}\right)+(1-a)\left[f_i\left(\frac{k+1}{3^i}\right)-f_i\left(\frac{k}{3^i}\right)\right],\\
&f_{i+1}\left(\frac{k+1}{3^i}\right)=f_i\left(\frac{k+1}{3^i}\right).
\end{flalign}
Given this construction, we say that
\begin{equation*}
F_a=\lim_{i\to\infty}f_i. 
\end{equation*}

We also can construct the graph of any $F_a$ with contraction mappings, using as a basis Katsuura's construction \cite{Katsuura91} of the graph of $F_{2/3}$.  Define contraction mappings $w_n: X \mapsto X$, where $n \in \{1,2,3\}$ and $X=[0,1] \times [0,1]$, as follows: For all $(x,y)\in X$,
\begin{flalign}
w_1(x,y)&=\left(\frac{x}{3},ay\right)\\
w_2(x,y)&=\left(\frac{2-x}{3},(2a-1)y+(1-a)\right)\\
w_3(x,y)&=\left(\frac{2+x}{3},ay+(1-a)\right).
\end{flalign}
In this case, the graph of $F_a$ is the unique invariant set for the iterated function system defined by (5), (6), and (7).

This paper builds primarily on a paper by Okamoto \cite{Okamoto05} regarding the same family of functions.  That paper showed, among other things, that if $a_0$ is the unique root of $27a^2-54a^3=-1$ in $(1/2,2/3)$, then $F_a$ is differentiable almost everwhere when $a\in(0,1/3)\cup(1/3,a_0)$; $F_{1/3}(x)\equiv x$ for $x\in[0,1]$ and thus is differentiable everywhere on that interval; $F_a$ is differentiable almost nowhere when $a\in(a_0,2/3)$; and $F_a$ is differentiable absolutely nowhere when $a\in[2/3,1)$. (Throughout this paper, ``almost everywhere'' and ``almost nowhere'' will refer specifically to sets of Lebesgues measure 1 and 0, respectively.)  Certain issues, however, were left unresolved in Okamoto's paper.  One such issue was whether $F_{a_0}$ is differentiable almost everywhere, almost nowhere, or neither.  Kobayashi \cite{Kobayashi09} recently proved that $F_{a_0}$ is in fact differentiable almost nowhere using the law of the iterated logarithm.

In this paper, we will give another proof of Kobayshi's result using only the properties of $F_a$ and the measure of the set of normal numbers in $[0,1]$.  Furthermore, using the same methods, we will describe how $F_a'$ behaves almost everywhere in $[0,1]$ as $a$ varies.  From here, we will prove that $F_a''$ does not exist anywhere for any $a\in(0,1)$ except $a=1/3$ (in which case $F_a(x)\equiv x$, and consequently, $F_a''(x)=0$ for all $x$) and $a=1/2$ (in which case $F_a'(x)=0$ for infinitely many disjoint subintervals of $[0,1]$).  Finally, we will show that if $a\in(0,1/2]$, then its graph has Hausdorff dimension 1, whereas if $a\in(1/2,1)$, then the Hausdorff dimension of its graph is equal to $\log_3 (12a-3)$.

\section{Remarks on differentiability}
\newtheorem{theorem}{Theorem}
\begin{theorem}
If $S_a$ is the set of all points $x$ for which $F_a'(x)$ does not exist and $a_0$ is the unique solution of $27a^2-54a^3=-1$ in $(1/2,2/3)$, then the Lebesgue measure $|S_{a_0}|=1$.
\begin{proof}
Let $S_a$ and $a_0$ be defined as they are above.  Now, we can write the ternary expansion of any $x\in[0,1]$ as
\begin{equation*}
x=0.\xi_1\xi_2\ldots=\sum_{i=1}^{\infty}\frac{\xi_i}{3^i}
\end{equation*}
where each $\xi_i\in\{0,1,2\}$. (If an expansion terminates, we say that it ends with a infinite string of zeroes. If $x=1$, we say that $\xi_i=2$ for all $i$.) For the first $n$ digits of the ternary expansion of $x$, we define $i(n)$ as the number of times $\xi_i=1$ in these first $n$ digits. Finally, we let $\displaystyle\gamma=\lim_{n\to\infty}\inf\frac{i(n)}{n}$.

According to Okamoto \cite{Okamoto05}, 
\begin{flalign*}
F_a'(x)&=\lim_{n\to\infty}(3-6a)^{i(n)}(3a)^{n-i(n)}\\
&=\lim_{n\to\infty}[3(1-2a)^\gamma a^{1-\gamma}]^n
\end{flalign*}
and thus $F_a'(x)=0$ if $|3(1-2a)^\gamma a^{1-\gamma}|<1$.  Setting $\gamma=1/3$, we can rewrite the second expression as $-1<3(1-2a)^{1/3} a^{2/3}<1$.  Raising each side to the third power, we get $F_a'(x)=0$ if $-1<27a^2-54a^3<1$; this quantity is equal to 1 when $a=1/3$, and it is equal to $-1$ when $a=a_0$.

As a result, when $a=a_0$ and $\gamma=1/3$, the sequence $\{f_i'(x)\}_{i=0}^\infty$ visits $1$ and $-1$ infinitely often, and therefore does not converge. \footnote{It should be noted, however, that $F_{a_0}'(x)$ does not converge absolutely to $1$, as we might expect; $\{f_i'(x)\}_{i=0}^\infty$ does not consist solely of `$1$'s and `$-1$'s, and as Kobayashi \cite{Kobayashi09} has pointed out, it contains subsequences that diverge absolutely to infinity.} Thus, $F_{a_0}'(x)$ does not exist for all $x$ satisfying $\gamma=1/3$.  To complete the proof, we must show that this set of values for $x$ has Lebesgue measure 1.

To do this, we turn to the \emph{normal numbers}---a set of irrational real numbers whose digits are uniformly distributed, regardless of the base in which the numbers are written.  Given any base $b\geq2$ and an alphabet of digits $\{0,1,\dots,b-1\}$, a normal number $x$ must satisfy
\begin{equation*}
\lim_{n\to\infty}\frac{i(n)}{n}=\frac{1}{b}
\end{equation*}
where $n$ represents the same value it did for $\gamma$ and $i(n)$ can be taken as the number of any individual digit in the first $n$ digits of $x$. (While this is not a sufficient property of normal numbers, it is a necessary one, and for the purposes of this proof, we only need to know that all normal numbers have this property.)

This means that the set of all $x\in[0,1]$ that satisfy $\gamma=1/3$ includes all the normal numbers in $[0,1]$.  Borel \cite{Borel09} proved that the set of all non-normal numbers has Lebesgue measure 0, so the set of all normal numbers in $[0,1]$ must have Lebesgue measure 1.  Because the normal numbers in $[0,1]$ are a subset of the numbers satisfying $\gamma=1/3$ in $[0,1]$, the latter set must have Lebesgue measure 1.  So $S_{a_0}$, the set of all $x\in[0,1]$ for which $F_{a_0}'(x)$ does not exist, must have Lebesgue measure $|S_{a_0}|=1$, as well.
\end{proof}
\end{theorem}

The fact that the normal numbers have Lebesgue measure 1 in $[0,1]$ reveals even more about $F_a'$ for all $a\in(0,1)$.  Because all normal numbers satisfy $\gamma=1/3$ in base 3, we know that the value of $F_a'(x)$ will depend on whether $|27a^2-54a^3|$ is less than, greater than, or equal to 1 for all $x$ in a set with Lebesgue measure 1.  This information allows us to give a general description of the behavior of $F_a'$ almost everywhere:
\newtheorem*{corollary}{Corollary}
\begin{corollary}
Let $a\in(0,1)$.
\begin{itemize}
\item If $a\in[a_0,1)$, then $F_a'(x)$ diverges for almost all $x\in[0,1]$;\\
\item If $a\in(0,1/3)\cup(1/3,a_0)$, then $F_a'(x)=0$ for almost all $x\in[0,1]$; and\\
\item If $a=1/3$, then $F_a'(x)=1$ for all $x\in[0,1]$.
\end{itemize}
\end{corollary}
We obtain the last result because $F_{1/3}(x)\equiv x$ and thus has the same derivative everywhere it is defined.  Additionally, if we take into account the fact that $F_a(x)$ is absolutely nowhere differentiable in $[0,1]$ when $a\in[2/3,1)$, we can replace \emph{almost all} with \emph{all} above when $a\in[2/3,1)$.

This leads us to our last question on the differentiability of $F_a$: What can we say about the second derivative $F_a''(x)$ for different values of $a$? We know that since $F_{1/3}(x)\equiv x$, $F_{1/3}''(x)=0$  for all $x$.  We also know that since $F_{1/2}'(x)=0$ for all $x$ except when $x$ belongs to the Cantor Set---a set with measure 0---$F_{1/2}'(x)$ is continuous on a set of points with measure 1.  Thus, $F_{1/2}'(x)$ is differentiable on this set of points, and its derivative $F_{1/2}''(x)=0$ at all of them.  But clearly $F_a''(x)$ does not exist at any $x$ when $a\in[a_0,1)$, since even $F_a'(x)$ fails to converge for almost all $x$ in this case.  The following theorem answers our question for the remaining values of $a$.

\begin{theorem}
For all $a\in(0,1/3)\cup(1/3,1/2)\cup(1/2,a_0)$ and all $x\in[0,1]$, $F_a''(x)$ does not exist.
\begin{proof}
Let $x\in[0,1]$.  We recall from Theorem 1 that almost all $x$ satisfy $\gamma=1/3$, in that the set of all such values has Lebesgue measure 1.  We also recall that for all $x$ belonging to this set of numbers and for all $a\in(0,1/3)\cup(1/3,a_0)$, $F_a'(x)=0$ if $|27a^2-54a^3|<1$. So for all $a\in(0,1/3)\cup(1/3,a_0)$ and almost all $x\in[0,1]$, $F_a'(x)=0$.

The graph of $y=F_a'(x)$ appears to be the line $y=0$ over $[0,1]$, but there are infinitely many points at which $F_a'(x)$ is undefined.  Okamoto \cite{Okamoto05} has shown that for $a\in(0,1/3)$, this set consists of all points of the form $x=(2k+1)/(2\cdot3^i)$, for $i\geq0$ and $k\in\{0,1,\dots,3^i-1\}$ and for $a\in(1/3,1/2)\cup(1/2,a_0)$, it consists of all points of the form $x=k/3^i$.  In either case, the set of points where $F_a'(x)$ does not exist is dense in $[0,1]$.  So $F_a'(x)$ is discontinuous everywhere in $[0,1]$, and therefore its derivative, $F_a''(x)$, does not exist anywhere in that interval.
\end{proof}
\end{theorem}

So unless $a=1/3$ or $a=1/2$---in which case $F_a''(x)=0$ for all $x$ at which it exists---the graph of $F_a$ is neither concave up nor concave down anywhere.  Also, $F_a$ belongs to differentiability class $C^0$ for all $a\neq1/3$.

Due to its somewhat pathological nature in terms of differentiability, this family of functions may call to mind another set of counterintuitive functions---that of the everywhere differentiable but nowhere monotone functions.  It is clear that these two sets have no common element; indeed, the only everywhere differentiable function in our set, $F_{1/3}$, is strictly increasing.  Nevertheless, some similarities exist in the ways these sets of functions behave on countable, dense sets.  For instance, if $a\in(0,1/3)$, then on the set of all points of the form $x=k/3^i$, where $i\geq0$ and $k\in\{0,1,\dots,3^i\}$, $F_a'(x)=0$, and at all points of the form $x=(2k+1)/(2\cdot3^i)$, $F_a'(x)$ is undefined; similarly, an everywhere differentiable but nowhere monotone function $G$ can be constructed over $[0,1]$ using the same two disjoint dense sets under the conditions that for all $x$ in one set, $G'(x)<0$ and that for all $x$ in the other set, $G'(x)>0$.  Beyond this, however, the author has found no connections between the constructions for these sets of functions.  This is primarily because with the dubious exceptions of $F_{1/3}$ and $F_{1/2}$, Okamoto's construction assures that if $F_a'(x)=0$ on one set dense in $[0,1]$, then $F_a'(x)$ must be undefined on another dense set.

\section{Remarks on dimension}
We first will show that because the graph of $F_a$ has finite arc length when $a\in(0,1/2]$, it has Hausdorff dimension 1.

\begin{theorem}
If $\Gamma_a$ is the graph of $F_a$, then its Hausdorff dimension $\dim_H(\Gamma_a)=1$ for all $a\in(0,1/2]$.
\begin{proof}
Suppose $\Gamma_a$ is the graph of $F_a$, and suppose $a\in(0,1/2]$.  Okamoto \cite{Okamoto05} has proven that for all $a\leq1/2$, $F_a$ is nondecreasing.  Now, if $F_a$ is nondecreasing, then clearly the affine pieces of each of its iteration $f_i$ must be nondecreasing as well, since $f_i(k/3^i)=F_a(k/3^i)$ for all $i\geq0$ and all $k\in\{0,1,\dots,3^i\}$.

Using the Triangle Inequality, we can determine the maximum arc length of $\Gamma_a$.  We do this by applying the inequality to each affine segment of $f_i$ with respect to the segment's horizontal and vertical components and letting $i\to\infty$.  More precisely, if $l_k$ is the length of the affine segment of $f_i$ over $[k/3^i,(k+1)/3^i]$, then $l_k\leq|f_i([k+1]/3^i)-f_i(k/3^i)|+|(k+1)/3^i-k/3^i|$.  

Because this inequality holds for every affine segment of $f_i$, we can obtain a bound for the arc length $\displaystyle L_i=\sum_{k=0}^{3^i-1}l_k$ of $f_i$.  And by taking the limit as $i\to\infty$, we can express the maximum arc length $L$ of $\Gamma_a$ as
\begin{equation*}
L\leq\lim_{i\to\infty}\sum_{k=0}^{3^i-1}\left|f_i\left(\frac{k+1}{3^i}\right)-f_i\left(\frac{k}{3^i}\right)\right|+\sum_{k=0}^{3^i-1}\left|\frac{k+1}{3^i}-\frac{k}{3^i}\right|.
\end{equation*}

Obviously, $(k+1)/3^i-k/3^i=1/3^i$ always, so the second sum becomes 1.  As for the first sum, we recall that the affine pieces of every $f_i$ are nondecreasing, so $|f_i([k+1]/3^i)-f_i(k/3^i)|=f_i([k+1]/3^i)-f_i(k/3^i)$.  Thus, for all $i$,
\begin{flalign*}
\sum_{k=0}^{3^i-1}\left[f_i\left(\frac{k+1}{3^i}\right)-f_i\left(\frac{k}{3^i}\right)\right]&=f_i\left(\frac{1}{3^i}\right)-f_i\left(\frac{0}{3^i}\right)+f_i\left(\frac{2}{3^i}\right)-f_i\left(\frac{1}{3^i}\right)+\cdots\\
&\hspace{1 pc}+f_i\left(\frac{3^i}{3^i}\right)-f_i\left(\frac{3^i-1}{3^i}\right)\\
&=f_i(1)-f_i(0)\\
&=1.
\end{flalign*}
So for $a\in(0,1/2]$, the maximum arc length of $\Gamma_a$ is 2.  To find a lower bound for the arc length, we first note that every $\Gamma_a$ has its endpoints at (0,0) and (1,1).  Because the shortest possible distance between these two points is a straight line, the minimum arc length of any $\Gamma_a$ must be $\sqrt{2}$.

Hence, the graph's 1-dimensional Hausdorff measure $\mathcal{H}^1(\Gamma_a)$ satisfies $0<\mathcal{H}^1(\Gamma_a)<\infty$, which implies that $\dim_H(\Gamma_a)=1$ for all $a\in(0,1/2]$.
\end{proof}
\end{theorem}

For $a\in(1/2,1)$, on the other hand, $F_a(x)$ is not nondecreasing, so we cannot find an upper bound for $L$ using the Triangle Inequality.  In fact, no such upper bound exists when $a\in(1/2,1)$; to prove this, we must show that the graph's Hausdorff dimension is greater than 1, but in order to do this, we first must calculate the box-counting dimension of the graph.

\begin{theorem}
If $\Gamma_a$ is the graph of $F_a$, then its box-counting dimension $\dim_B(\Gamma_a)=\log_3(12a-3)$ for all $a\in(1/2,1)$. 
\begin{proof}
Suppose $\Gamma_a$ is the graph of $F_a$, and suppose $a\in(1/2,1)$.

To determine the box-counting dimension of $\Gamma_a$, we will use a variant of the more familiar method of box-counting.  Instead of counting squares of side length $\delta$, we will determine the smallest area $A(f_i)$ needed to cover the graph of iteration $f_i$ using rectangles that share a horizontal side length of $\delta$.  Then the ``number'' of squares with side length $\delta$ needed to cover the graph of $f_i$ can be expressed as $N_\delta(f_i)=A(f_i)/\delta^2$.  This variant works by the same principle as box-counting dimension with squares does, but it gives more precise answers for each $\delta$ we use.

Letting $\delta=1/3^i$, we see that for all $\Gamma_a$, $A(f_0)=1$.  Using the horizontal and vertical ratios of contraction for (5), (6), and (7), we see that
\begin{flalign*}
A(f_{i+1})&=\frac{1}{3}aA(f_i)+\frac{1}{3}(2a-1)A(f_i)+\frac{1}{3}aA(f_i)\\
&=A(f_i)\left(\frac{4a-1}{3}\right).
\end{flalign*}
So by induction,
\begin{equation*}
A(f_i)=\left(\frac{4a-1}{3}\right)^i
\end{equation*}
and thus,
\begin{flalign*}
N_\delta(f_i)&=\frac{([4a-1]/3)^i}{(1/3)^{2i}}\\
&=(12a-3)^i.
\end{flalign*}
So for $a\in(1/2,1)$, the box-counting dimension of $\Gamma_a$ is
\begin{flalign*}
\dim_B(\Gamma_a)&=\lim_{i\to\infty}\frac{\log([12a-3]^i)}{-\log(1/3^i)}\\
&=\log_3 (12a-3).
\end{flalign*}
\end{proof}
\end{theorem}

Next, we will show that the Hausdorff dimension of $\Gamma_a$ is equal to its box-counting dimension by the Mass Distribution Principle (see, e.g. \cite{Falconer03}):

\begin{theorem}
If $\Gamma_a$ is the graph of $F_a$, then its Hausdorff dimension $\dim_H(\Gamma_a)=\log_3(12a-3)$ for all $a\in(1/2,1)$.
\begin{proof}
We consider a small variation on Katsuura's construction of $\Gamma_a$ using contraction mappings.  Let $E_0=[0,1]\times[0,1]$ and define further levels of the construction by $E_{i+1}=w_1(E_i)\cup w_2(E_i)\cup w_3(E_i)$ where $i>0$ and $w_1$, $w_2$, and $w_3$ are mappings (5)--(7).  Clearly $E_{i+1}\subset E_i$ for all $i\geq0$, and $\displaystyle\bigcap_{i=0}^\infty E_i=\Gamma_a$.  This last relationship between the $E_i$s can be understood as follows: While in Okamoto's construction, linear segments are constructed ``upwards'' to a graph with infinite length, in this construction, rectangular regions are constructed ``downwards'' to the same graph which has zero area.  Moreover, we see that each $E_i$ can be covered by $3^i$ rectangles of length $(1/3)^i$.

Using methods related to the box-counting process in Theorem 4, it can be shown that if $a\in(1/2,1)$, then the area of $E_{i+1}$ can be expressed as
\begin{flalign*}
A(E_{i+1})&=a\left(\frac{1}{3}\right)A(E_i)+(2a-1)\left(\frac{1}{3}\right)A(E_i)+a\left(\frac{1}{3}\right)A(E_i)\\
&=\frac{4a-1}{3}A(E_i),
\end{flalign*}
and since $A(E_0)=1$, we have $\displaystyle A(E_i)=\left(\frac{4a-1}{3}\right)^i$ for all $i\geq0$.

Now, let $\mu$ be the natural mass distribution on $\Gamma_a$; we start with unit mass on $E_0$  and repeatedly ``spread'' this mass over the total area of each $E_i$.  Also, let $U$ be any set whose diameter $|U|<1$.  Then there exists some $i\geq0$ such that
\begin{equation*}
\left(\frac{1}{3}\right)^{i+1}\leq|U|<\left(\frac{1}{3}\right)^i,
\end{equation*}
an inequality that applies to any $U$ satisfying $0<|U|<1$.  Given these conditions on the diameter of $U$, it is clear that for every $U$, there is some $i$ such that $U$ is contained in an open square of side length $(1/3)^i$ and $U$ contains points in at most two level-$i$ ``sub-rectangles.''

Hence, the area of $U$ is bounded above by the area of the open square containing it; that is, $A(U)\leq(1/9)^i$.  In terms of measure, we know that the entire area of $U$ can be contained in $E_i$, so 
\begin{flalign*}
\mu(U)&\leq\frac{A(U\cap E_i)}{A(E_i)}\\
&\leq\frac{(1/9)^i}{([4a-1]/3)^i}\\
&\leq\left(\frac{1}{12a-3}\right)^i.
\end{flalign*}
And since $\displaystyle\left(\frac{1}{3}\right)^{i+1}\leq|U|$ implies that $\displaystyle\left(\frac{1}{3}\right)^i\leq3|U|$, we have
\begin{equation*}
\mu(U)\leq\left(\frac{1}{12a-3}\right)^i=\left(\frac{1}{3^i}\right)^{\log_3 (12a-3)}\leq(3|U|)^{\log_3 (12a-3)}=(12a-3)|U|^{\log_3 (12a-3)},
\end{equation*}
and therefore, by the Mass Distribution Principle, $\log_3 (12a-3)\leq\dim_H(\Gamma_a)\leq\dim_B(\Gamma_a)$, and given the upper bound obtained in Theorem 4, we have $\dim_H(\Gamma_a)=\log_3 (12a-3)$.
\end{proof}
\end{theorem}

This means that if $a>1/2$, then $\dim_H(\Gamma_a)>1$.  And if the Hausdorff dimension of $\Gamma_a$ exceeds 1, then its one-dimensional Hausdorff measure---that is, its length---must be infinite.  So for $a>1/2$, the graph's arc length $L$ indeed has no upper bound.

Given the result of Theorem 5, it is clear that $\dim_H(\Gamma_a)$ ranges continuously between 1 and 2 for all $a\in(1/2,1)$.  This means that for $a\in(1/2,1)$, $\Gamma_a$  meets the criterion of a fractal, being a self-similar set whose Hausdorff dimension exceeds its topological dimension.

\section{Conclusions}
In his paper, Okamoto \cite{Okamoto05} posed another question related to $F_a$'s differentiability: If $a\in(0,1/3)$, there are infinitely many $x$ for which $F_a'(x)$ is finite but nonzero, but how do we classify them?  In the corollary to Theorem 1, we have provided a partial answer to this question.  Since $F_a'(x)=0$ if $a\in(0,1/3)$ and $x$ satisfies $\gamma=1/3$, we know that the points in question at least must be non-normal numbers in $[0,1]$.  It follows that the set of such points must occupy a subset of $[0,1]$ with Lebesgue measure 0.

\section{Acknowledgments}
The author would like to thank Professor Daniel Jackson for direction and assistance with this project.


\begin{thebibliography}{25}
\bibitem{Borel09}
{\scshape Borel, \'{E}mile.}
Les probabilit\'{e}s d\'{e}nombrables et leurs applications arithm\'{e}tiques.
\emph{Rendiconti del Circolo Matematico di Palermo} \textbf{27.1} (1909) 247--71.
\bibitem{Bourbaki04}
{\scshape Bourbaki, Nicolas.}
Functions of a real variable: elementary theory.
Trans. from the 1976 French original [MR 0580296] by Philip Spain. 
\emph{Springer-Verlag, Berlin}, 2004. MR 2013000.
\bibitem{Falconer03}
{\scshape Falconer, Kenneth.}
Fractal geometry---mathematical foundations and applications. Second edition. 
\emph{John Wiley \& Sons, Inc., Hoboken, NJ}, 2003. MR 2118797 (2006b:28001).
\bibitem{Katsuura91}
{\scshape Katsuura, Hidefumi.}
Continuous nowhere-differentiable functions---an application of contraction mappings.
\emph{American Mathematical Monthly} \textbf{98.5} (1991) 411--416. MR 1104304 (92c:26005).
\bibitem{Kobayashi09}
{\scshape Kobayashi, Kenta.}
On the critical case of Okamoto's continuous non-differentiable functions.
\emph{Proc. Japan Acad. Ser. A Math. Sci.} \textbf{85.8} (2009) 101--104. MR 2561897 (2010m:26006).
\bibitem{McCollum}
{\scshape McCollum, James.}
Properties of Bourbaki's function.
\bibitem{Okamoto05}
{\scshape Okamoto, Hisashi.}
A remark on continuous, nowhere differentiable functions.
\emph{Proc. Japan Acad. Ser. A Math. Sci.} \textbf{81.3} (2005) 47--50. MR 2128931 (2006a:26005), Zbl 1083.26004.
\end{thebibliography}
\end{document}